\documentclass[11pt]{amsart}
\usepackage[left=1.25in,right=1.25in,top=1in]{geometry}
\usepackage{graphicx,verbatim,array}
\usepackage{graphics,url}
\usepackage[normalem]{ulem}
\usepackage{amssymb}
\usepackage{amsmath,amsthm}
\usepackage{multirow} 
\usepackage{etoolbox,booktabs}

\usepackage{amsmath, amsthm}

\newtheorem{theorem}{Theorem}
\newtheorem{corollary}[theorem]{Corollary}
\newtheorem{lemma}{Lemma}
\newtheorem{proposition}{Proposition}

\theoremstyle{remark}
\newtheorem*{rem}{Remark}

\theoremstyle{definition}

\usepackage{amsmath, amsthm}
\usepackage{color}

\newcommand{\R}{\mathbb{R}}

\def\z1{z_{1}}

\def\Diff{\mathbf{Diff}}

\newcommand{\mathd}{\mathrm{d}}
\newcommand{\dx}[1]{\mathd #1}

\newcommand{\ppx}[1]{\frac{\partial}{\partial #1}}

\newcommand{\pfpx}[2]{\frac{\partial #1}{\partial #2}}

\usepackage[textwidth=1.25in,textsize=tiny]{todonotes}
\usepackage[footnotesize,bf]{caption}

\usepackage{siunitx}
\title[Computation of Weil-Peterson Geodesics]{Numerical Computation of Weil-Peterson Geodesics in the Universal
Teichm\"uller Space}
\author[Feiszli]{Matt Feiszli}
\author[Narayan]{Akil Narayan}
\thanks{A. Narayan was supported by National Science Foundation Awards DMS-1318427 and DMS-1552238}

\renewcommand{\ker}{\mathrm{ker}}

\newcommand{\pinv}{\dagger}

\newcommand{\UTS}{PSL_{2}(\R) \backslash \Diff(S^{1})}

\newcommand{\lmpt}{\mathcal{Q}}

\newcommand{\lm}{\mathcal{L}^{M}}

\newcommand{\V}[1]{\mathbf{#1}}

\usepackage{ulem}
\usepackage{algorithm,algorithmic}

\begin{document}
\begin{abstract}
We propose an optimization algorithm for computing geodesics on the universal Teichm\"uller space $T(1)$ in the Weil-Petersson ($W P$) metric. Another realization for $T(1)$ is the space of planar shapes, modulo translation and scale, and thus our algorithm addresses a fundamental problem in computer vision: compute the distance between two given shapes. The identification of smooth shapes with elements on $T(1)$ allows us to represent a shape as a diffeomorphism on $S^1$. Then given two diffeomorphisms on $S^1$ (i.e., two shapes we want connect with a flow), we formulate a discretized $W P$ energy and the resulting problem is a boundary-value minimization problem. We numerically solve this problem, providing several examples of geodesic flow on the space of shapes, and verifying mathematical properties of $T(1)$. Our algorithm is more general than the application here in the sense that it can be used to compute geodesics on any other Riemannian manifold.
\end{abstract}
\maketitle

\section{Introduction and Background}

Representation and comparison of shapes is a central problem in computer vision.  In the past several decades, many approaches to represent, compare, and classify shapes have been presented (See 
\cite{veltkamp_shape_2001}, 
\cite{mumford_pattern_2002} for review and discussion).  The space of 2D shapes is inherently nonlinear; this poses fundamental difficulties in computer vision when attempting object recognition and statistics.  In \cite{sharon_2d-shape_2006}, Mumford and Sharon describe a construction based on conformal mapping which makes the space of simple closed plane curves into a Riemannian metric space.  The space itself is in fact the universal Teichm\"uller space with the Weil-Peterson metric; in this paper, we describe a numerical solver for geodesics in this space. In Sections \ref{sec:intro-welding} and \ref{sec:intro-wp} we describe the optimization problem we wish to solve: minimization of an energy functional with boundary value constraints. We note that a minimization algorithm proposed in \cite{sharon_2d-shape_2006} was applied to only relatively simple shapes because of numerical difficulties. In this work we aim to apply our algorithm to more general, complicated shapes. Our work is competitive with a recent approach based on shooting for the analogous boundary value problem \cite{kushnarev-2012}.

The outline of this paper is as follows: Section \ref{sec:discrete} describes the mathematics of computations for the Weil-Peterson metric using discrete samples of a velocity field. Section \ref{sec:computations} discusses the fully discrete algorithm along with our technique for satisfying the boundary constraints. Finally, Section \ref{sec:results} presents numerical results that illustrate the effectiveness of the algorithm.


\subsection{Conformal Welding}\label{sec:intro-welding}

The construction in \cite{sharon_2d-shape_2006} is classical in Teichm\"uller theory.  Given a simply-connected planar region $\Omega$ bounded by a smooth Jordan curve (this is our definition of a ``shape''), one constructs a pair of conformal maps $\Phi_{+}: \Delta \to \Omega_{+}$ and $\Phi_{-}: \Delta \to \Omega_{-}$, which map the exterior and interior of the unit disk to the exterior and interior of $\Omega$, respectively.  The exterior map is normalized to fix $\infty$ and have real derivative there.  The interior map is only defined up to right multiplication by the three-parameter M\"obius group $PSL_{2}(\R)$ of conformal self-maps of the unit disk.  Both $\Phi_{+},\Phi_{-}$ extend continuously to the boundary $S^{1}$, and their composition $\Phi_{-}^{-1} \circ \Phi_{+}$, restricted to the boundary, is a map $\Psi: S^{1} \to S^{1}$.

\begin{figure}
  \begin{center}
    \resizebox{1.0\textwidth}{!}{
    \includegraphics{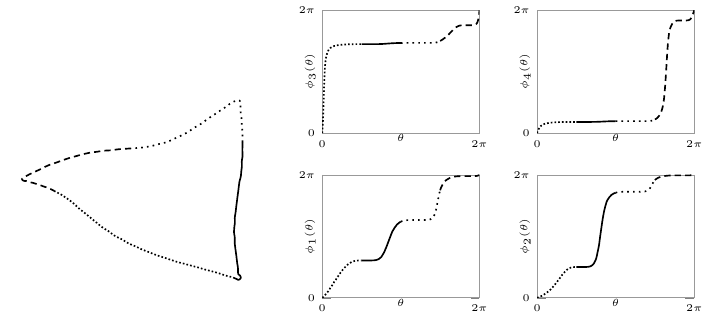}
    }
  \end{center}
  \caption{A triangular region $\Omega$ (left) transformed into an equivalence class of welding maps $\phi$ (right). Four members of the equivalence class are shown. Four segments of $\partial \Omega$ are drawn with differing linestyles to highlight their corresponding segment on each of the welding maps.}
  \label{fig:welding-maps} 
\end{figure}

A remarkable result (see 
\cite{ahlfors_lectures_2006}, for example) is that this result is almost an isomorphism: the space of shapes, modulo translation and scale, is isomorphic to the group $\Diff(S^{1})$, modulo conformal self-maps of the disk.  This provides an elegant way of making the space of shapes into a metric space:  we take an element of the coset space $\UTS$, i.e. an equivalence class of diffeomorphisms of the circle, as a representation of a shape.  This space is known as the universal Teichm\"uller space and was initially studied in the context of Riemann surfaces 
\cite{ahlfors_lectures_2006,hubbard_teichmuller_2006}
it also arises  in string theory 
\cite{bowick_string_1987,nag_rm_1990}.  In \cite{sharon_2d-shape_2006} elements of $\UTS$ are called ``fingerprints''; in the mathematical literature, elements of the broader class of quasisymmetric homeomorphisms of $S^{1}$ are commonly known as ``welding maps'' and each has an associated quasicircle.  We illustrate an example of a equivalence class of welding maps in Figure \ref{fig:welding-maps} for a simple shape.

\subsection{The Weil-Peterson Riemannian metric}\label{sec:intro-wp}

The Lie algebra to $\UTS$ is the space of vector fields on the circle modulo the subspace spanned by $1, \cos$ and $\sin$.
If $v(\theta) = \sum_{n=-\infty}^{\infty} a_{n} e^{in\theta} \frac{\partial}{\partial \theta}$ is a real-valued vector field on the circle, the Weil-Peterson norm is 
\begin{align}\label{eq:wp-norm-fourier}
\| v \|_{WP}^{2} = \sum_{n=2}^{\infty} (n^{3}-n) |a_{n}|^{2} = \frac{1}{2} \left( \left\| v^{(3/2)} \right\|^2 - \left\|v^{(1/2)}\right\|^2 \right)
\end{align}
where $\left\|\cdot\right\|$ denotes the $L^2$ norm.  \footnote{This formula, while explicit, is troublesome in numerical computations; the boundary values of conformal maps have so much high-frequency content, even for curves with real-analytic boundaries, that the number of Fourier coefficients required for for an accurate global representation is prohibitive.}  This may be rewritten as
\begin{align*}
\| v \|_{WP}^{2} = \langle Lv, v \rangle
\end{align*}
for the positive semidefinite, self-adjoint operator
\begin{align}\label{eq:operator-L-definition}
L = -\mathcal{H}(\partial^{3} - \partial)
\end{align}
where $\mathcal{H}$ is the Hilbert transform.  Note that $L$ has a kernel which is exactly the span of $1, \cos$ and $\sin$.  As mentioned above, the corresponding vector fields in the span of $\left\{ \sin\theta \frac{\partial}{\partial \theta}, \cos\theta \frac{\partial}{\partial \theta}, \frac{\partial}{\partial \theta}  \right\}$  are infinitesimal M\"obius maps; i.e. they span the Lie algebra $psl_{2}(\R)$ of the M\"obius group $PSL_{2}(\R)$.  Hence, our norm is indeed a norm on the quotient space.

By right-translations we extend this norm to the entire space; that is, if $\phi(\theta, t) \equiv \phi_{t}(\theta)$ is a curve in $\UTS$, then any any time $t$ we pull back the velocity $\dot \phi_{t}$ under the derivative of right translation to obtain
\begin{align*}
v_{t}(\theta) = \dot \phi_{t} \circ \phi_{t}^{-1}
\end{align*}
The length of a curve $\phi:[0,1] \times S^{1} \to S^{1}$ in $\Diff(S^{1})$ is obtained by integrating
\begin{align*}
L[\phi] = \int_{0}^{1} \| \dot \phi_{t} \circ \phi_{t}^{-1} \|_{WP} dt
\end{align*}
The resulting Riemannian metric is known as the Weil-Peterson metric.  The classical definition was in terms of Beltrami differentials of quasiconformal self-maps of the hyperbolic plane and was originally studied in the context of Teichm\"uller theory.  The formulation we present here occurs naturally in string theory where it was discovered as the unique K\"ahler metric on $\UTS$ 
\cite{segal_definition_2004}.  See 
\cite{nag_rm_1990} for a thorough discussion of how this K\"ahler metric agrees with the classical WP metric.

It is a fact from variational calculus that minimizing length is equivalent to minimizing energy: the path of least energy is a constant-speed parametrization of the path of least length.  Thus given a pair of points $h_{0}, h_{1} \in \Diff(S^{1})$ we find the geodesic by finding the path $\phi$ of minimum energy
\begin{align}\label{eq:energy-definition}
E[\phi] = \int_{0}^{1} \| \dot \phi_{s} \circ \phi_{s}^{-1} \|_{WP}^{2} ds,
\hskip 10pt \phi_0 = h_0, \hskip 10pt \phi_1 = h_1
\end{align}
Finding this minimizing $\phi$ is the problem we are concerned with.  It is a fact that all sectional curvatures of the WP metric are negative 
\cite{teo_weil-petersson_2006}, which in finite dimensions would imply that minimizing geodesics are unique.  While the existence and uniqueness of minimizing geodesics is a subtle issue in infinite dimensions, both questions were recently answered in the affirmative in 
\cite{gay-balmaz_geometry_2009}:  Any two points in $\UTS$ are joined by a unique geodesic, and further, solutions to the geodesic equation exist for all time and the resulting welding maps are Sobolev $H^{s}$ for any $s < 3/2$, where the inequality is sharp.

It is an essential fact that the WP metric we have described here is, by construction, \textit{right-invariant}.  That is, given any $h_{0}, h_{1}$, define the map $g$ by $h_{0} = g \circ h_{1}$ and let $\phi(\theta, s)$ be the shortest path from $g$ to the identity.  Right-invariance of the metric means that the length of $\phi(\theta, s)$ is the same as length of $\phi(h_{1}(\theta), s)$, which immediately implies that the shortest path from $h_{0}$ to $h_{1}$ is in fact $\phi(h_{1}(\theta), s)$.  Hence, computing geodesics between two welding maps reduces to the case where one of the maps is the identity.   Therefore in what follows we consider only paths ending at the identity.

\subsection{Path energy and its gradient}\label{sec:pathlength-notation-intro}

In order to compute the gradient of path energy, we first need to put a metric on the space of admissible paths.  We can then take the gradient in this metric (as opposed to simply taking the gradient in our particular choice of coordinate).  Choosing a good metric greatly improves convergence and stability.  This is sometimes called the ``natural'' gradient, particularly in the machine learning community. 

Let $\mathcal{P}$ be the space of smooth paths on $\UTS$.  Let $\phi(s,\theta, t) = \phi_{s}(\theta, t)$ be a smooth curve in $\mathcal{P}$.  Variable $\theta$ is position on $S^{1}$, $t$ parametrizes a path for given $s$, and $s$ parametrizes the curve of paths.  Then for any $s$ the path of a particle $q(s, \theta, \cdot)$ is obtained by integrating the Eulerian vector field $v_{s}(\theta, t)$ defined by:
\begin{align*}
\phi_{s}(\theta, t) = \phi_{s}(\theta, 0) + \int_{0}^{t} v_{s}(\phi_{s}(\theta, \xi), \xi) d\xi
\end{align*}
Denote the variation as
\begin{align}\label{eq:pathspace-velocity-intro}
w_{s}(\theta, t) = \ppx{s}\left( \dot{\phi_{s}} \circ \phi_{s}^{-1}\right)(\theta, t)
\end{align}
and introduce the norm
\begin{align}\label{eq:pathspace-norm-intro}
\left\| w_{s} \right\|_{\mathcal{P}}^{2} \equiv \int_{0}^{1} \| w_{s} \|_{WP}^{2} dt
\end{align}
where the WP norm is taken in $\theta$.  It is in this metric on path space that we shall take the gradient. 

We make one further remark:  there is a natural identification between paths $\phi$ which end at the identity and their velocity fields $v_{t} = \dot{\phi_{t}} \circ \phi_{t}^{-1}$.  In what follows, we prefer to take velocity fields $v(t, \theta)$ as coordinates on path space, as opposed to working with paths themselves.  We will then compute the natural gradient of energy for our discrete boundary-value problem in two steps: we first compute the gradient for the unconstrained problem; i.e. we compute an update to the velocity field which simply makes the energy smaller, ignoring the boundary conditions at times $0$ and $1$.  We then project this update (in the metric induced by \eqref{eq:pathspace-norm-intro}) onto the space of admissible updates; i.e. those updates which preserve the boundary conditions.  This two-step trick is equivalent to simply taking the first variation of energy for the fixed-endpoint problem, but is considerably simpler to derive, carries  little performance penalty, and, while experimenting, eliminates the requirement to recalculate the constrained gradient when the method for computing the unconstrained gradient is changed.  The significant difficulty of performing numerical work in quotient spaces made this last point quite practical.

\subsection{Example}

At this point it might be helpful to illustrate our method with a simplified example.  Since all sectional curvatures of the WP metric are negative, we make our example on the hyperbolic plane: a manifold of constant negative curvature. Since geodesics in this model are explicitly known and computable, this serves as a straightforward test case for the algorithm. We emphasize that in this example we make simple choices for the discretization in order to make the presentation clear. Finally we note that this example does not showcase some of the more nuanced aspects of our algorithm. (E.g., quotient spaces, more sophisticated quadrature, complicated boundary conditions, etc.) These more subtle points will be discussed later in the article when our algorithm is applied to the WP metric.

\begin{figure}
  \includegraphics[width=\textwidth]{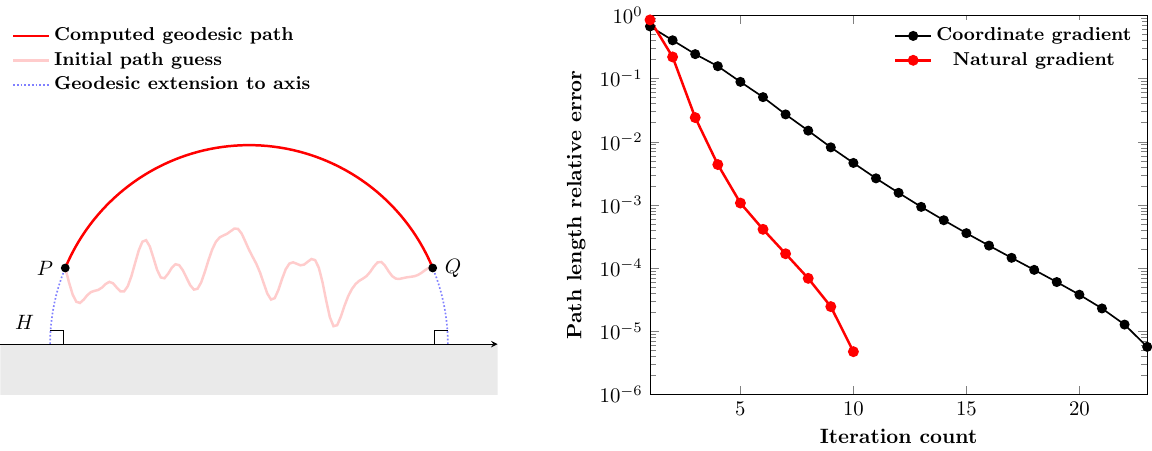}
  \caption{Left: Computed geodesic path in the hyperbolic plane model, along with a randomly-generated initial guess. Right: Path length (versus exact value) for the natural-gradient and coordinate-gradient algorithms, as a function of iteration count. The noise level is $\tau = 0.05$.}\label{fig:hyperbolic-plane}
\end{figure}

In the half-plane model of the hyperbolic plane, the set of points
\begin{align*}
H = \{ (x,y) \in \R^2 \; |\; y > 0 \}
\end{align*}
is equipped with the distance element
\begin{align}\label{eq:half-plane-ds}
ds^{2} = \frac{dx^{2} + dy^{2}}{y^{2}}
\end{align}
Geodesics in this model are circular arcs orthogonal to the x-axis; see Figure \ref{fig:hyperbolic-plane}, left. 

Consider two points $P$ and $Q$ in $H$, and let $\left( (x_0, y_0),\, (x_1, y_1),\, \ldots, (x_N, y_N) \right)$ denote an ordered set that form any $(N+1)$-point initial guess to a geodesic connecting $P$ and $Q$. We constrain the endpoints:
\begin{align*}
  P &= (x_0, y_0), & Q &= (x_N, y_N).
\end{align*}
For now we assume that the points along the initial path $q_n = (x_n, y_n)$ for $n=1, \ldots, N-1$ are initialized in any way. A basic discretization involves a velocity field collocated at the midpoint of $q_{n-1}$ and $q_n$:
\begin{align}\label{eq:half-plane-v}
  v_{n-\frac{1}{2}} &\triangleq q_n - q_{n-1} \in \R^2, & n &= 1, \ldots, N
\end{align}
Our algorithm actually works on the velocity $v_{n-\frac{1}{2}}$ and uses it to define and update the path $q_n$\footnote{Note that the transformation $q_n \rightarrow v_{n-\frac{1}{2}}$ corresponds to a size $N-1$ linear and invertible transformation: There are $N$ total velocity variables $v_{n-\frac{1}{2}}$, and $N+1$ total position variables $q_n$. Adding the boundary conditions constrains 1 velocity variable and 2 position variables.}. This discretization implies that the path length energy is approximated by summing up discrete versions of \eqref{eq:half-plane-ds} collocated at midpoints between the $q_n$:
\begin{align}\label{eq:half-plane-E}
  E \approx E_N &= \sum_{n=1}^N \dx{s}_{n-\frac{1}{2}}^2 = \sum_{n=1}^N \frac{\|v_{n-\frac{1}{2}}\|^2}{ y^2_{n - \frac{1}{2}} },  \\\nonumber
  y_{n - \frac{1}{2}} &= \frac{1}{2} \left( y_{n+1} + y_n \right)
\end{align}
Since $y_n$ (i.e., $q_n$) is a function of $v_{n-\frac{1}{2}}$, this implies that the equation above is implicitly a function of $v_{n-\frac{1}{2}}$. Therefore, we can take the gradient of $E_N$ in the $v_{n-\frac{1}{2}}$ coordinates $\nabla E_N = \left( \pfpx{E_n}{v_{n-\frac{1}{2}}} \right)_{n=1}^N$. Energy is minimized along geodesics, so we can perform numerical gradient descent with the $v_{n-\frac{1}{2}}$ as variables to minimize $E_N$. 

However, coordinate-dependent gradients, like the one considered above, are not the natural direction of descent for optimization algorithms in non-Euclidean geometries.  When the manifold in question has a known metric (or a suitable metric can be imposed), it is preferable to take the gradient in this metric, and this is sometimes referred to as the \textit{natural} gradient \cite{amari_natural_1998}. In our example, the distance element is proportional to $1 / y$, so updating the $\{ q_{n} \}$ according the to coordinate gradient permits particles near the boundary to take arbitrarily large steps.  In this case, \eqref{eq:half-plane-E} implies that in our $v_{n-\frac{1}{2}}$ coordinate system, a diagonal matrix with entries $G_{n,n} = y^{-2}_{n - \frac{1}{2}}$ will give us a metric tensor on path space which will cause path updates to respect the metric. Therefore, the direction $G^{-1} \nabla E_N$ is the steepest descent direction.  

However, this natural gradient direction does not respect the endpoint conditions. (I.e., updating $v_{n-\frac{1}{2}}$ and then computing the path positions $q_n$ via \eqref{eq:half-plane-v} will not in general satisfy both $P = q_0$ and $Q = q_N$. To correct for this, we project $G^{-1} \nabla E_N$ into the space of admissible updates that do satisfy the endpoint constraints. We omit presenting the explicit details of this here, but the general continuous version of the projection is \eqref{eq:operator-Y-definition} discussed in Section \ref{sec:pathlength-notation}. The discrete (matrix) version of this projection in the WP metric case is given in \eqref{eq:discrete-path-metric}.

The high-level description of the algorithm can now be completed: at each stage with a given $v_n$ we compute the constrained natural gradient of $E_N$, and use it to perform steepest descent; we iterate until convergence.

For our numerical examples here, we initialize the geodesic path guess as a random perturbation of a Euclidean straight line:
\begin{align*}
  q_n &= P + \frac{n}{N+1} \left(Q - P\right) + \sigma_n, & n &= 1, \ldots, N-1
\end{align*}
where $\left(\sigma_n\right)_{n=1}^{N-1}$ is treated as a Gaussian process with covariance structure $\mathrm{cov} \left( \sigma_n, \sigma_m \right) = \tau \exp\left( -\frac{(n-m)^2}{(N + 1)\ell} \right)$, with correlation length $\ell$ and noise magnitude $\tau$. In our examples below, we set $\ell = \frac{1}{500}$, and $\tau$ varies between $0$ and $0.1$. (See Figure \ref{fig:hyperbolic-plane} left, for an example realization of this process.) 

To illustrate the effectiveness of the natural gradient versus the coordinate gradient we set $\tau = 0.05$ and $N=100$, and compare the different gradient iteration schemes in Figure \ref{fig:hyperbolic-plane} right. The natural gradient converges to the true geodesic much more quickly than the coordinate gradient.  In fact, in flat spaces it will converge in a single iteration.

\begin{figure}
  \includegraphics[width=\textwidth]{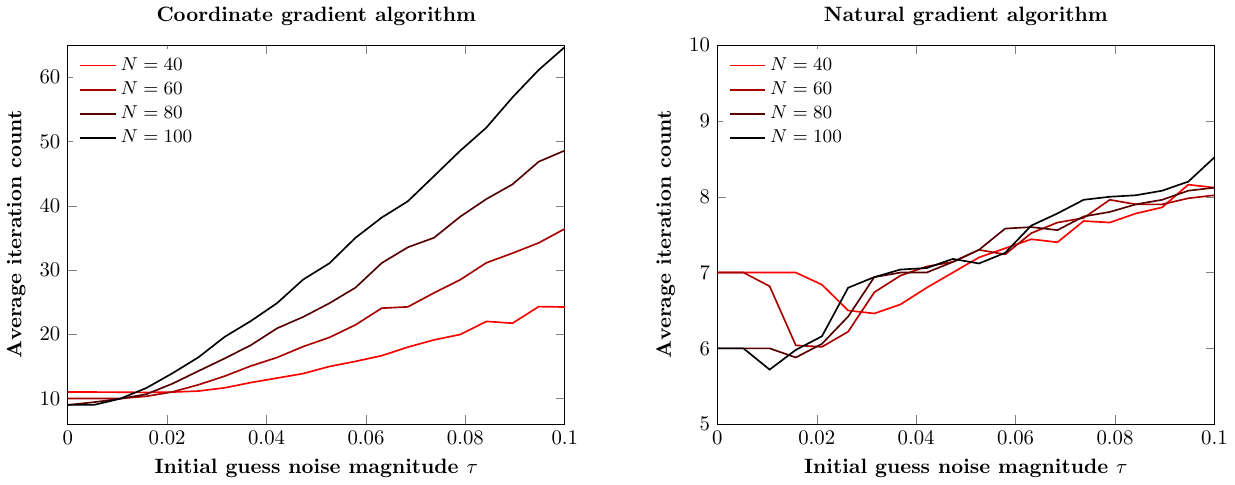}
  \caption{Ensemble-averaged iteration counts required to converge to hyperbolic plane geodesic from a randomly-perturbed Euclidean straight line initial guess (see Figure \ref{fig:hyperbolic-plane}, left). Left: Coordinate gradient algorithm. Right: Natural gradient algorithm.}\label{fig:hyperbolic-plane-mc}
\end{figure}

The natural gradient is also more robust with respect to discretization error and initial guesses. In Figure \ref{fig:hyperbolic-plane-mc} we show the required number of iterations until convergence for various values of $N$ and $\tau$. Since our initial guess is a random process, we run a size-$50$ ensemble of simulations, and the plots show ensemble average iteration counts. The natural gradient algorithm is quite insensitive to both noise level and discretization parameters, whereas the coordinate gradient is quite sensitive and also requires more iterations.  In what follows we will use a very similar metric to stabilize updates to our discrete path in shape space.

\section{The Discrete Problem}\label{sec:discrete}


We now return to the problem \eqref{eq:energy-definition} in the case where $h_{0}: S^{1} \to S^{1}$ is an arbitrary welding map and $h_{1}$ is  the identity map.  This is a boundary-value minimization problem: find $v_s = \dot{\phi}_s \circ \phi_s^{-1}$ subject to initial and
terminal conditions $\phi_0 = h_{0}$ and $\phi_1 = \mathrm{Id}$ such that $\int \|v_s
\|_{WP}^2 \dx{s}$ is minimized. We introduce $q \equiv \phi$, and given a velocity field $v$ then 
\begin{align*}
  q(t,x) = q(0,x) + \int_0^t v(\xi, q(\xi,x)) \dx{\xi}
\end{align*}
for all $0 \leq t \leq 1$.

We track $M$ particles on the range of $q$ as they evolve in $t$.
For some fixed $x_m$ on $S^1$, $q(t, x_m)$ are the locations of these particles. The velocity field $v(t, \cdot)$ at the locations $q(t,x_m)$ is likewise stored.
We store these velocity field values and recover the particle positions $q(t,x_m)$ by integrating rather than directly storing the particle positions. 
In the next subsection we will describe how we compute the $WP$ norm at each time. We delay discretizing the $t$ variable until Section \ref{sec:computations}.


\subsection{Computing the WP Norm}

We first consider the problem of computing the norm of a vector field on $S^1$.  Subsequent sections will consider computation of path energy and its gradient to find geodesics.   

Our discretization provides us with information about the velocity field only at a discrete set of points.  Our first task is to interpolate this velocity field to a field defined on all of $S^{1}$ and then compute the norm of this lift. The computation proceeds in 3 steps.  We first collect some standard results which let us do this when the metric has no kernel.  Next we compute the optimal interpolant and its norm when the metric has a kernel; this is done by projecting the interpolation data onto the orthogonal complement of the kernel and computing the norm of the lift.  Finally, we generalize the problem to allow an arbitrary basis for the interpolating functions; this allows for more flexibility in the implementation and also suggests a more stable method of computation, which we detail later.

\subsubsection{Lifting the velocity field into the orthogonal complement of the kernel}\label{sec:lifting-v}

We need to extend a vector field defined on a finite subset of $S^{1}$ to a vector field on all of $S^{1}$.  Here we describe a method for finding the smoothest possible extension using Green's functions for the operator $L$.  In this section we consider only vector fields orthogonal to the kernel; we extend to more general cases next.

Consider the space of all vector fields on $S^1$, and a subspace on which
$\|\cdot\|_{W P}$ is a proper norm:
\begin{align*}
  V &= \{v: \|v \|_{W P} < \infty\}, & \widetilde{V} = \{v \in V: v \perp 
  \ker_V L\},
\end{align*}
where, explicitly, $\ker_V L = \mathrm{span}\{1, \sin, \cos\}$.
Let $\lm$ be the $M$-manifold of configurations of $M$ points $q_{1} <
q_{2} < ... < q_{M}$ on the unit circle and write $\lmpt = \left\{
q_{1}, q_{2}, ..., q_{M}  \right\}$. The tangent space
$T_{\lmpt}\lm$ is the $M$-dimensional space of vector fields $v = (v_{1},
v_{2}, ..., v_{M})$ supported on the set $\lmpt$.  We consider all possible
extensions of $v$ to vector fields defined on the entire circle, and induce a
norm on $T_{\lmpt}\lm$ by 
\begin{align}\label{eq:v-global-minimizer}
\| v \|_{WP(\lmpt)} \equiv \inf \left\{ \| \widetilde v \|_{WP} : \widetilde{v}(q_{m}) = v_{m}, 1 \le m \le M  \right\}
\end{align}
Conversely, given the basepoint $\lmpt$, consider the evaluation map taking a global vector field $\widetilde v$, defined on all of $S^{1}$, to $T_{\lmpt}\lm$ by
\begin{align*}
\widetilde v \to \left( \widetilde v(q_{1}), \widetilde v(q_{2}), ..., \widetilde v(q_{M})   \right)
\end{align*}
This map induces a splitting of any $\widetilde v$ into two components.  The set of vector fields for which $\widetilde v(q_{m}) = 0$ for all $1 \le m \le M$ is the \textit{vertical subspace} (of the \textit{Lie algebra} of $\Diff(S^{1})$: the set of smooth vector fields on $S^{1}$), and its $WP$-orthogonal complement is the \textit{horizontal subspace}; the minimizing extension in \eqref{eq:v-global-minimizer} of $v \in T_{\lmpt}\lm$ is known as a \textit{horizontal lift}.

As the next proposition shows, the horizontal subspace is an $M$-dimensional subspace spanned by translates of Green's function for the operator $L = -\mathcal{H}(\partial^{3} - \partial)$.  This Green's function is known explicitly, see \cite{kushnarev_teichons_2009}:
\begin{align*}
G(\theta) = 2 \sum_{n=2}^{\infty} \frac{\cos(n\theta)}{(n^{3}-n)} = (1-\cos \theta) \log\left[ 2(1-\cos \theta)  \right] + \frac{3}{2}\cos \theta -1
\end{align*}
The following results are standard facts which may be proven using the reproducing kernel property of Green's function and the fact that the vector fields $G(\theta-q_{m})\frac{\partial}{\partial\theta}$ lie in the horizontal subspace.
\begin{proposition}[Computing Horizontal lifts]
\label{prop:horizontal-lift-nokernel}
\ \newline
\begin{enumerate}
\item The set of vector fields
\begin{align*}
\left\{ G(\theta - q_{m}) \frac{\partial}{\partial\theta} \right\}_{m=1}^{M}
\end{align*}
is a basis for the component of the
horizontal subspace in $\widetilde{V}$.
\item The horizontal lift on $\widetilde{V}$ of the tangent vector $(v_{1}, ..., v_{M})
\in T_{\lmpt}\lm$ is the vector field $\widetilde v(\theta)
\frac{\partial}{\partial \theta}$ for
\begin{align*}
\widetilde v(\theta) = \sum_{i=1}^{M}G(\theta - q_{i})p_{i} \frac{}{}\quad where \quad p_{i} = \sum_{j=1}^{M}G^{-1}_{ij}v_{j}
\end{align*}
and $G_{ij}$ is the positive definite symmetric matrix
\begin{align}\label{eq:G-matrix}
G_{ij} = \left\langle G(\cdot - q_i), G(\cdot - q_j) \right\rangle_{W P} = G(q_{i} - q_{j} )
\end{align}
Further,
\begin{align*}
\| \widetilde v \|_{WP}^{2} = \sum_{i,j}G^{-1}_{ij} v_{i}v_{j}
\end{align*}
\end{enumerate}
\end{proposition}

\subsubsection{Lifting the velocity field}
Now we extend these results to compute horizontal lifts which take into account the kernel.

\begin{lemma}\label{lemma:v-minimizer}
Let $G$ be an $M\times M$ positive matrix, $B$ an $M \times K$ full-rank matrix
with $K < M$. Then given some $v \in \R^M$, the projection of $v$ onto the
$G^{-1}$-orthogonal complement of the kernel of $B^T$ is given by 
\begin{align*}
  P_G^B v \doteq (I - B (B^T G^{-1} B)^{-1} B^TG^{-1}) v
\end{align*}
\end{lemma}
In our case $B$ is an $M \times 3$ matrix containing point evaluations of a basis for the three-dimensinal kernel of the $W P$ norm.  Once we have subtracted out the contribution of the kernel, we may compute the lift and its norm.  The next result is a corollary of section \ref{sec:lifting-v} and the previous lemma.

\begin{figure}
  \begin{center}
    \resizebox{1.0\textwidth}{!}{
      \includegraphics{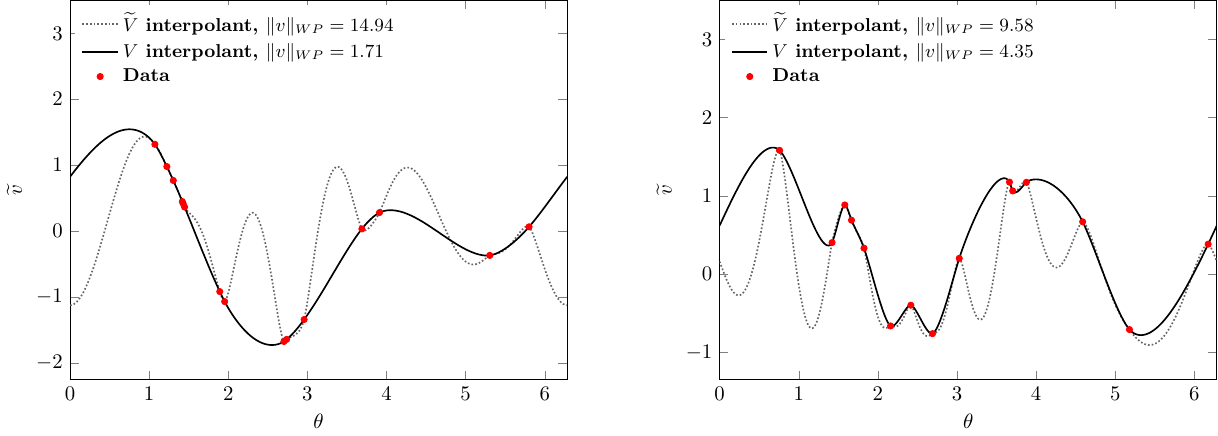}
    }
  \end{center}
  \caption{Interpolation procedures from (a) Proposition \ref{prop:horizontal-lift-nokernel} on the quotient space $\widetilde{V}$, and from (b) Corollary \ref{cor:v-field-wp-min} on $V$. Left: the data is taken from the function $v(\theta) = \sin (2 \theta) + \cos(\theta)$ on 15 randomly distributed $\theta$ locations. Right: the data is taken as Gaussian random variable perturbations of function values from $v(\theta) = \sin (2 \theta)$ on 15 randomly distributed $\theta$ locations.}
  \label{fig:minimal-interpolation}
\end{figure}

\begin{corollary}\label{cor:v-field-wp-min}
Fix $\lmpt$ and let the matrix $G$ be as in Proposition
\ref{prop:horizontal-lift-nokernel}. For any $v \in T_\lmpt \lm$, 
\begin{align}\label{eq:v-discrete-minimizer}
  \| v \|^2_{WP(\lmpt)} = v^{T}G^{-1}P_{G}^{B}v  
\end{align}
Further, the horizontal lift $\widetilde{v}$ of $v$ to $V$ is given by
\begin{align*}
  \widetilde{v} &= \sum_{j=1}^3 w_j
  b_j(\cdot)  + \sum_{j=1}^M \left( G^{-1} P_G^Bv
  \right)_j G(\cdot - q_j)
\end{align*}
where
\begin{align*}
  w &= (B^T G^{-1} B)^{-1} B^T G^{-1}v 
\end{align*}
and $\{b_j\}_{j=1}^3$ is the basis for $\ker_V L$ used to construct $B$:
$B_{m,k} = b_k(q_m)$. 
\end{corollary}
\begin{rem}
The vector $p = G^{-1}P_{G}^{B}v$ contains what are known as the momenta coefficients. Similarly,we call $w$ the kernel coefficients.
\end{rem}
The interpolant produced from Corollary \ref{cor:v-field-wp-min} is notably different than the horizontal lift from the quotient space in Proposition \ref{prop:horizontal-lift-nokernel}. This is exemplified in Figure \ref{fig:minimal-interpolation}.

Up to this point, we have taken Green's functions centered at the interpolation nodes as 
as the basis for our interpolating space.  As discussed, this produces the norm-minimal lift.  While this optimality is nice, and is in fact what our current implementation ultimately uses, one might wish to choose a different space of interpolating functions.  For example, the first $N$ complex exponentials, or a wavelet basis adapted to the problem might be reasonable choices.  We address this now; it is only a small modification to the calculations above and admits a stable solution which we present later.

Let $\widetilde{V}_N$ be any $N$-dimensional subspace of $\widetilde{V}$ such that the interpolation problem for data collocated at $\lmpt$ is unisolvent.  If
$F = \left\{ f_n(\cdot) \right\}$ is a basis of $\widetilde{V}_N$, define the $M \times N$
matrix $\lambda_F$ to have entries 
\begin{align*}
\left(\lambda_F\right)_{m,n} = \langle f_n, G(\cdot-q_m)\rangle_{W P} = f_n(q_m)
\end{align*}
Let $G_F$ be the $N \times N$ Gram matrix for the $f_{n}$:
 \begin{align*}
\left(G_F\right)_{m,n} = \langle f_m, f_n\rangle_{WP}
\end{align*}
 For $v \in \R^M$ interpolation data on $\lmpt$, define $\widetilde{V}_N(v)$ as
the subset of functions from $\widetilde{V}_N \oplus \ker_V L$ that interpolate
to $v$. Then define
\begin{align*}
  \| v \|^2_{WP(\lmpt), \widetilde{V}_N} = 
  \min_{v \in \widetilde{V}_N \oplus \ker_V L} \| v\|^2_{W P}
\end{align*}
The results from the previous section may be rephrased by replacing $G \gets G_{F}$
where $G_{F} = \lambda_F G_F^{-1} \lambda_F^T$.  Since $\widetilde{V}_N$ is total for interpolation, then $G_F$ is invertible. One finds that all formulas carry over essentially up to change-of-basis.  This, then, is our most general result for computing lifts, which we state in terms of an operator $L_{F}$.

\begin{proposition}\label{prop:metric-definition}
In the notation above, set
\begin{align*}
L_{F} = G_{F}^{-1} P^B_{G_{F}}
\end{align*}
Then 
\begin{align}\label{eq:other-subspace-norm}
  \| v \|_{WP(\lmpt), \widetilde{V}_N}^2 
                &= v^T L_{F} v \\\nonumber
                &=: v^T p.
\end{align}
The reconstruction $\widetilde{v} \in \widetilde{V}_N \oplus \ker_V
L$ is given by 
\begin{align*}
  \widetilde{v} = \sum_{n=1}^N c_n f_n + \sum_{j=1}^3 w_j b_j
\end{align*}
where
\begin{subequations}
\label{eq:other-subspace-coefficients}
\begin{align}
  w &= \left(B^T G_{F}^{-1}  B \right)^{-1} B^T G_{F}^{-1} v \\
  c &= G_F^{-1} \lambda_F^T G_{F}^{-1} (v - B w).
\end{align}
\end{subequations}
\end{proposition}

Note that unlike the result of Corollary \ref{cor:v-field-wp-min}, the momenta
 $p \equiv L_{F}v$ are no longer the reconstruction coefficients for the lift
 unless $\widetilde{V}_N$ is the span of Green's functions centered at $\lmpt$. In any case, $L_F$ is a symmetric, rank $M-3$, positive semidefinite matrix.

Whether using Green's function interpolation or another basis, using
\eqref{eq:other-subspace-norm} and \eqref{eq:other-subspace-coefficients} to
compute either the norm or any coefficients can be problematic. The matrices in
the explicit formula given can easily be ill-conditioned, depending on $F$ and
$\mathcal{Q}$. We do not directly solve this issue here, but instead we reformulate the problem in
a format that we observe is more robust with respect to numerical precision issues than the above formulas. This is done in Section \ref{sec:computations}.

\subsection{Path energy and its gradient}\label{sec:pathlength-notation}

Computing the length and energy of a path is relatively straightforward given the work in the previous section. If $v^t$ is a vector of velocity evaluations at time $g$, we have
\begin{align*}
  E = \int_0^1 \|\widetilde{v}\|^2_{W P} \dx{t} \int_0^1 v^t L_F^t v^t \dx{t},
\end{align*}
where $L^t_F$ is the metric on the particle positions $q^t$.

As mentioned above, we compute the gradient of energy for our discrete boundary-value problem in two steps.  We first compute the gradient ignoring the boundary conditions at times $0$ and $1$.  We then project this update onto the manifold of updates which preserve the boundary conditions. Since we have yet to discretize the temporal variable, we wait to compute the uncontrained gradient until Section \ref{sec:computations}. For now we concentrate on the projection onto the admissible updates.

Define the energy of a velocity field as above.  Now consider two manifolds of
velocity fields.  $\mathcal{N}$ contains all velocity fields defined on $M$
particles at each time, and $\mathcal{M} \subset \mathcal{N}$ contains those
which respect the boundary-value conditions.
\begin{align*}
\mathcal{N} &= \left\{ \widetilde{v}:\{1,...,M\} \times [0,1] \to \R     \right\} \\
\mathcal{M} &= \left\{ V \in \mathcal{N} \left| q_m(1) = q_m(0) + \int_0^1 \widetilde{v}(\xi,q_m(\xi)) \dx{\xi} \right. \right\}
\end{align*}

We defined the energy $E: \mathcal{N} \to \R$ as a function on $\mathcal{N}$, and in the next section we will describe the computation of its gradient $\nabla E$ on $\mathcal{N}$.  
With this gradient in hand, we obtain the direction of steepest ascent as a projection of $\nabla E$ onto the manifold $\mathcal{M}$ by projection in the following metric.

Let $\mathcal{P}$ be the space of smooth paths on $\UTS$.  Let $q(t,x,s) = q_{t;s}(x)$ be (a representative of) a smooth curve in $\mathcal{P}$.  Variable $x$ is position on $S^{1}$, $t$ parametrizes a path for given $s$, and $s$ parametrizes the curve of paths.  Recall that the velocity $v_{t;s}$ along the path is defined by
\begin{align*}
q(t,x,s) = q(0,x,s) + \int_{0}^{t} v_{\xi;s}(q(\xi,x,s)) d\xi
\end{align*}
and as described above, we write
\begin{align}\label{eq:pathspace-velocity}
w_{t;s} = \ppx{s}\left( \dot{q} \circ q^{-1}\right)_{t;s}
\end{align}
and introduce the norm
\begin{align}\label{eq:pathspace-norm}
\left\| w_{;s} \right\|_{\mathcal{P}}^{2} \equiv \int_{0}^{1} \| w_{t;s} \|_{WP}^{2} dt
\end{align}

However, we wish to work with the velocity field manifolds $\mathcal{N}$ and $\mathcal{M}$ as our coordinate on $\mathcal{P}$.  With $q = q(t,x,s)$, $\dot{q} = \ppx{t} q$, and $q' = \ppx{x} q$, we have 
 \begin{align*}
\pfpx{\left(\dot{q} \circ q^{-1}\right)_{t;r}}{r}(x) &= \left( \pfpx{\dot{q}}{s} \circ q^{-1}\right)_{t;r} + (\dot{q} \circ q^{-1})'_{t;r} \left(\pfpx{q^{-1}}{s}\right)_{t;r} \\
&= \left(\pfpx{v}{s}\right)_{t;r} + v_{t;r}' \left(\pfpx{q^{-1}}{s}\right)_{t;r} \\
&= \left(\pfpx{v}{s}\right)_{t;r} - v_{t;r}' \left(q^{-1}\right)'_{t;r} \left(\pfpx{q}{s}\circ q^{-1}\right)
 \end{align*}
 Direct computation yields
 \begin{align*}
 \frac{\partial q_{s}}{\partial s}(t,x) &= \frac{\partial}{\partial s} q(0,x,s) + \frac{\partial}{\partial s}\int_{0}^{t} v_{\xi;s}(q(\xi,x,s)) d\xi \\
 &= \int_{0}^{t} \frac{\partial}{\partial s} v_{\xi;s}(q(\xi,x,s)) d\xi 
 \end{align*}
 Therefore
 \begin{align*}
\pfpx{\left(\dot{q} \circ q^{-1}\right)_{t;r}}{r}(x) &= \mathcal{Y}\left[ \pfpx{v}{s}\right],
 \end{align*}
 where the linear operator $\mathcal{Y}$ is given by 
 \begin{align}\label{eq:operator-Y-definition}
\mathcal{Y}(t,x,s) = \mathrm{Id} - v' \left(q^{-1}\right)' \int_0^t (\cdot)_{\xi} \dx{\xi}
 \end{align}

This operator $\mathcal{Y}$ gives a mapping from velocity field updates (i.e. tangent vectors to $\mathcal{N}$ or $\mathcal{M}$) to tangent vectors $\frac{\partial q}{\partial s}$ to $\mathcal{P}$.  Hence our metric on path space immediately induces a metric on $\mathcal{N}$ and the submanifold $\mathcal{M}$.  We make our projection of the update to the velocity field with respect to this metric.

\section{Computations}\label{sec:computations}
We now address the more practical issues of implementation of the methods described. We first formulate computation of the $W P$ norm as a classical numerical linear algebra problem. Following this is presentation of the temporal discretization; having discretized both the temporal and $S^1$ variables, we can compute the unconstrained coordinate gradient of the energy. Finally, the submanifold projection for gradient updates is discretized.

\subsection{The WP norm}
It is not necessary to explicitly form the metric matrix $L_F$ as defined in Proposition \ref{prop:metric-definition}. We can accomplish the same task in a more stable manner using some established numerical linear algebra results. Given an $M
\times N$ matrix $A$, $A^\pinv$ denotes its Moore-Penrose pseudoinverse. For any
$\V{b} \in \R^N$, the vector $A^\pinv \V{b}$ is the minimum-norm, least squares
solution to $A \V{x} = \V{b}$.  Here, both minimum-norm and least
squares refer to the $\ell^2$ (Euclidean) metric on $\V{x}$ and $\V{b}$,
respectively. 

For computation of the $W P$ norm, we develop the ideas for a general subspace
$V_N$ of $V$ and then specialize to
the minimal-norm lift over all $V$ using Green's functions. Our goal is the interpolation of data
$\V{v} \in \R^M$ collocated at $\lmpt$ with an element from
$\widetilde{V}_N \oplus \ker_V L$; thsu we are essentially trying to find the minimal norm solution to
\begin{align}\label{eq:full-minnorm-system}
  \left[\begin{array}{cc} \lambda_F & B \end{array}\right] \left(
  \begin{array}{c} \V{c} \\ \V{w} \end{array} \right) =:
  A \V{x} = 
  \V{v},
\end{align}
where the norm over the coefficients $\V{c}$ is given by the positive (Gram)
matrix $G_F$ and that over the kernel coefficients $\V{w}$ is zero. Again we
assume that $\widetilde{V}_N$ is total for interpolation on $\lmpt$,
implying that $A \V{x} = \V{v}$ has at least one solution. Let $G_{F_0}$ be a
block diagonal matrix with $G_F$ in the upper-left block and a $3 \times 3$
zero matrix in the lower-right block.  It is easy to show that $\ker (B)$ is
trivial so long as $M \geq 3$.  This implies that $\ker(G_{F_0}) \cap \ker(A)$
is trivial and therefore the minimum norm solution to
\eqref{eq:full-minnorm-system} is unique \cite{eldn_weighted_1982}.  This is an
easy way to show uniqueness of the lift defined in
\eqref{eq:v-global-minimizer}.

Restating \eqref{eq:full-minnorm-system}, we are trying to solve the
minimization problem:
\begin{align}\label{eq:lse-problem}
  \text{minimize} \hspace{0.5cm} \left\| R \V{x} \right\|_2
  \hspace{0.5cm} \text{subject to} \hspace{0.5cm}
  A \V{x} = \V{v},
\end{align}
where $R$ is any matrix square root of $G_{F_0}$. (I.e.
$R^T R = G_{F_0}$.) One of the standard tools for solving
least squares equality constrained problems is the generalized singular value
decomposition (GSVD) \cite{loan_generalizingsingular_1976}. Since $A$ and $G(F)_0$ have the same number of
columns, the GSVD matrix decomposition of $(A, R)$ is given by 
\begin{subequations}
\label{eq:gsvd}
\begin{align}
  A &= U C X^{-1} \\
  R &= V S X^{-1},
\end{align}
\end{subequations}
where $U$ and $V$ are orthogonal matrices, $X$ is invertible, and $C$ and $S$
each have only one non-vanishing diagonal with non-negative entries and satisfy
$C^T C + S^T S = I$. (Here, $V$ denotes a matrix in the GSVD decomposition, and not a linear subspace.) The solution to \eqref{eq:lse-problem} is given by
\begin{align}\label{eq:gsvd-solution}
  \V{x} = \left(\begin{array}{c} \V{c} \\ \V{w}\end{array}\right) &= X C^\pinv U^T \V{v},
\end{align}
The above matrix mapping $\V{v}$ to $\V{x}$ is a weighted pseudoinverse of $A$
with respect to $R$.  The momenta are $\V{p} = U C^{\pinv T} X^T G_{F_0}
\V{x} = L_F \V{v}$.  Our experience is that this is a much more stable way to compute
momenta than the explicit matrix relations used to define $L_F$. Note that
since $C^\pinv$ is a matrix with only one non-vanishing diagonal, applications of
$C^{\pinv}$ can be accomplished with a simple vector-vector multiply. The norm is
$\|\V{v}\|_{WP(\lmpt),\widetilde{V}_N} = \V{c}^T G_F\V{c} = \V{v}^T \V{p}$.

For the horizontal lift from the whole space $V$, we need only take the space
$\widetilde{V}_N$ to be the $M$-dimensional space of Green's functions centered
at the $q_m$. Notable simplifications in this case are that (a) $G_F = G$,
with $G$ being the Gram matrix of Green's functions \eqref{eq:G-matrix} and (b) the momentum and basis
coefficients coincide: $\V{p} = \V{c}$.

\subsection{Temporal discretization}
Since flow along the geodesic is reversible, we use a quadrature method that is symmetric with respect to the endpoints. We rewrite the flow of particles in an equivalent form
\begin{align}\label{eq:symmetric-quadrature}
  q_m(t) = \frac{1}{2} \left(q_m(0) + \int_0^t \widetilde{v}(\xi, q_m(\xi)) \dx{\xi}\right) + \frac{1}{2} \left(q_m(1) - \int_t^1 \widetilde{v}(\xi, q_m(\xi)) \dx{\xi}\right)
\end{align}
Now choose $T$ ordered points $s^t$ on $(0,1)$; let $\V{q}^t$ and $\V{v}^t$ denote the $M$ particle positions and velocities, respectively, at those times. The endpoint particle positions $\V{q}^0$ and $\V{q}^{T+1}$ are given.
Let $h^t$ denote the quadrature weight associated with data at time $s^t$ so that $\int_0^1 f(s) \dx{s} \simeq \sum_{t=1}^T h^t f(s^t)$.

Following \eqref{eq:symmetric-quadrature} we use the quadrature rule to
to integrate particle positions forward from $\V{q}^0$ to $\V{q}^t$, and backwards from $\V{q}^{T+1}$ to $\V{q}^t$, and average the result. For example if the representation is piecewise-linear, then the symmetric velocity-to-particle map is 
\begin{align}\label{eq:explicit-particle-update}
  \V{q}^t &= \frac{1}{2}\left(\V{q}^0 + \sum_{r=1}^{t-1} h^r \V{v}^r +
  \frac{h^t}{2} \V{v}^t\right) + \frac{1}{2}\left( \V{q}^{T+1} - \sum_{r=t+1}^T
  h^r \V{v}^r - \frac{h^t}{2} \V{v}^t\right), & t\geq 1
\end{align}
The scheme described above expresses particle positions $\V{q}^t$ linearly with the velocity field
values $\V{v}^r$. Collect all the velocity
evaluations into a matrix $V$, of size $M \times T$; do the same for the
particle positions in a matrix $Q$.  Then regardless of the choice of linear temporal quadrature rule, there exists a $T \times T$ matrix $Z$
such that 
\begin{align}\label{eq:particle-update}
  Q = \frac{1}{2} \left(Q^0 + Q^{T+1}\right) + V Z,
\end{align}
where $Q^0$ and $Q^{T+1}$ are matrices with the known vectors $\V{q}^0$ and
$\V{q}^{T+1}$ repeated. The entries of the matrix $Z$ depends on the choice of temporal representation and quadrature. For the piecewise-linear choice, it has entries $Z_{r,t} = \frac{1}{2} c_{r,t} h_r$, where 
\begin{align}\label{eq:crt-def}
c_{r,t} = \left\{\begin{array}{rl} 
                   0, & r = t, \\
                   1, & r < t, \\
                   -1, & r > t.
                 \end{array}
          \right.
\end{align}
In order to evaluate the energy $\int_0^1 \| v_s\|_{W P}^2 \dx{s}$, we build the quadrature factor $h^t$ into the norm at each point in time. This can be accomplished by simply replacing $G(F)$ by $h^t G(F)$ at each point in time. Therefore
\begin{align}\label{eq:pathlength-discretization}
  E = \sum_{t=1}^T \| \V{v}^t \|_{WP(\lmpt),\widetilde{V}_N}^2 = \sum_{t=1}^T 
  (\V{v}^t)^T L_{F} \V{v}^t 
\end{align}
where $L_{F}$ is time-varying, depending explicitly on $\V{q}^t$,
the particle locations at time $s^t$, and also depending proportionally on $h^t$.

To keep the presentation simple we have described use of a piecewise-linear
quadrature. In practice we find it is sometimes more efficient to use a Legendre-Gauss-Lobatto quadrature rule for $t
\in [0,1]$, which enables a high-order polynomial representation for each
particle's velocity field. When to use a high-order representation versus a piecewise linear representation depends on whether one expects large $t$-derivatives in the velocity field.

\subsection{The gradient of the WP norm}
At some time $s^t$, we have seen that if $L_{F} = G_{F}^{-1}P_{G_{F}}^{B}$ then
formula \eqref{eq:pathlength-discretization} defines the norm, and to minimize
we must compute variations of this with respect to velocity evaluations. Since
the particle positions are influenced by \eqref{eq:particle-update}, the
matrices $L_F$ and therefore the energy change in nontrivial ways when we vary
$\V{v}^t$. To simplify the procedure, we first fix $t$ and consider only
variations in $\|\V{v}^t\|_{WP(\lmpt)}$ with respect to $\V{v}^r$ for $r =
1,\ldots,T$.

Direct variation of the quadratic form yields
\begin{align}\nonumber
  \ppx{\V{v}^r} \|\V{v}^t\|^2_{WP(\lmpt),\widetilde{V}_N} &= 2 \delta_{t,r} (\V{v}^t)^T L_{F} + (\V{v}^t)^T \left(\ppx{\V{v}^r} L_F \right) \V{v}^t \\ \label{eq:wpnorm-variation-velocity}
  &= 2 \delta_{t,r} (\V{p}^t)^T + (\V{v}^t)^T \left(\ppx{\V{v}^r} L_F \right) \V{v}^t
\end{align}
where $\delta_{i,j}$ is the Kronecker delta.
The metric $L_F$ and its components $G_{F}$ and $P_{G_{F}}^B$ are generated from $\V{q}^t$, depending on $\V{v}^t$. Using simple properties of
matrix algebra, a straightforward but messy computation yields a formula for the
variation of the objective with respect to the particle positions.
\begin{lemma}\label{lemma:wpnorm-particle-variation}
Let $\lambda(F')$ be the $M \times N$ matrix with entries
$\lambda(F')_{m,n} = f_n'(q_m)$ and $B'$ be the $M \times 3$ matrix
with entries $B'_{m,j} = b_j'(q_m)$.  For two matrices $C$ and $D$ of the
same size, $C \circ D$ denotes the elementwise product.  Then 
\begin{align}\label{eq:wpnorm-variation-particle}
 \ppx{\V{q}^t} \|\V{v}^t\|^2_{WP(\lmpt),\widetilde{V}_N} = -2 \left[ \V{p}^t \circ
 \left(\lambda(F') G(F)^{-1} \lambda(F)^T \V{p}^t + B'
 \V{w}^t\right)\right]^T.
\end{align}
In particular, if $F$ is a basis of Green's functions centered at $\lmpt$, then 
\begin{align}\label{eq:wpnorm-variation-particle-hlift}
 \ppx{\V{q}^t} \|\V{v}^t\|^2_{WP(\lmpt)} = -2 \left[ \V{p}^t \circ
 \widetilde{v}'(\V{q}^t)\right]^T,
\end{align}
where $\widetilde{v}'$ is the derivative of the horizontal lift from $V$.
\end{lemma}
The salient result of Lemma \ref{lemma:wpnorm-particle-variation} is that
computing variations of the $W P$ norm with respect to particles is quite easy
once we have the momentum coefficients $\V{p}$ and kernel coefficients $\V{w}$
in hand from the GSVD. The only additionaly difficulty could come in computing
entries for $\lambda(F')$. In all the straightforward choices for $F$ we make,
it is not a deterrent. 

\subsection{The gradient of energy}
With the gradient of the $W P$ norm computed, computing the energy is now
straightforward. Combining \eqref{eq:particle-update},
\eqref{eq:wpnorm-variation-velocity}, and \eqref{eq:wpnorm-variation-particle},
we have
\begin{align}\label{eq:standard-gradient-formula}
\frac{1}{2} \pfpx{E}{\V{v}^t} =  (\V{p}^t)^T - 
                      \sum_{r}  Z_{t,r} \left[ \V{p}^t \circ
 \left(\lambda(F') G(F)^{-1} \lambda(F)^T \V{p}^t + B'
 \V{w}^t\right)\right]^T
\end{align}
where again we recall that the
momenta coefficients $\V{p}^t$ already have the scaling factor $h^t$ included.

Although the above equation is an exact formula for the gradient, updating
velocity fields with these values will not respect the endpoint constraint
$\V{q}^1 = \V{q}^0 + \sum_t h_t \V{v}^t$, and so we project into the space of admissible updates.

The discrete versions of the manifolds presented in Section \ref{sec:pathlength-notation} can be defined as follows: 
\begin{align*}
N &= \left\{ \V{v}:[1,...,M] \times [1, ..., T] \to \R     \right\} \\
M &= \left\{ \V{v} \in N \left| \V{q}^1 = \V{q}^0 + \sum_t h_t \V{v}^t \right. \right\}
\end{align*}
Since the condition $\V{q}^1 = \V{q}^0 + \sum_t h_t \V{v}^t$ is linear, we may write it as
\begin{align*}
\V{q}^1 - \V{q}^0 = A \V{v}
\end{align*}
for some constraint matrix $A$. We see that $M$ is an implicit submanifold of
$N$: it is the $\V{q}^1 - \V{q}^0$ level set of the function $A \V{v}$.  By the
implicit function theorem, the tangent space to $M$ at any point is equal to
the nullspace of the differential $DA$, and we have simply $[DA] = A$. That is,
the admissible updates to a given velocity field $V$ are exactly those vectors
$\V{v}$ in the nullspace of $A$.

In order to perform the appropriate projection onto the nullspace of $A$, we place our update along a curve on the space of paths $\mathcal{P}$ and apply the norm induced by the metric \eqref{eq:pathspace-velocity} and \eqref{eq:pathspace-norm}. We then seek a discretization of the semidefinite operator $\mathcal{Y}^\dagger L \mathcal{Y}$ given by \eqref{eq:operator-Y-definition} and \eqref{eq:operator-L-definition}, where $\mathcal{Y}^\dagger$ is the $L^2$ adjoint of $\mathcal{Y}$.

At a fixed value of $s^t$, we have already approximated $L$ with $L_F$ depnding
on $\V{q}^t$.  Therefore we concentrate on $\mathcal{Y}$: the integral
$\int_0^{s^t} (\cdot)_\xi \dx{\xi}$ is approximated by the $t$-th column of
$Z$. Let $\widetilde{Z} = Z^T \otimes I_M$, where $I_M$ is the $M\times M$ identity
matrix; then $\V{q} = \widetilde{Z} \V{v}$. 
The factor $v'(\V{q})$ can be
computed at each $t$ via $v'(\V{q}^t) = \lambda(F') \V{c}^t + B' \V{w}^t$. We
collect these factors into a diagonal $M T \times M T$ matrix $W$ with the
$t$-th block diagonal entry $(W^t)_{j,j} = v'(q^t_j) (q^{-1})'(q^t_j)$.  The
result then is the modified metric on updates $\ppx{\V{v}} E$:
\begin{align}\label{eq:discrete-path-metric}
  \widetilde{L_F} = (I - W \widetilde{Z})^T L_F (I - W \widetilde{Z}).
\end{align}
I.e., the matrix $Y := I - W \widetilde{Z}$ is our approximation to the operator $\mathcal{Y}$.
We use this metric both to form the natural gradient and to orthogonally project into the nullspace of the constraint matrix $A$.

\section{Algorithm and Results}\label{sec:results}
With all the necessary derivations complete, a summary of the minimization algorithm is presented in Algorithm Listing \ref{alg:min-algorithm}.
One detail we have not yet addressed is a method for computing welding maps -- i.e. for determining $\phi$ given a discrete collection of ordered samples on $\partial \Omega$. A few standard methods exist for numerically computing conformal welds \cite{driscoll_schwarz-christoffel_2002,li_numerical_1998,symm_integral_1966,symm_numerical_1967}, and we settle upon the Zipper algorithm \cite{marshall_convergence_2007}; specifically we use the (simplest) `geodesic' version. We use the Zipper algorithm to (a) construct a representation of a conformal weld given discrete samples on the shape boundary $\partial \Omega$, (b) `interpolate' the welding map to any starting particle locations $\V{q}^0$,$\V{q}^{T+1}$ of our choosing, and (c) `invert' the weld: compute a representation of a simple closed curve from discrete samples of a homeomorphism of $S^1$. There are other possibilities for accomplishing these conformal welding tasks \cite{sharon_2d-shape_2006}.

\begin{algorithm}
\begin{algorithmic}

\STATE Input: Initial and final particle positions $\V{q}^0$ and $\V{q}^{T+1}$
\STATE Initialize: $\V{v}^t = \V{q}^{T+1} - \V{q}^0$ for all $t$
\WHILE{Not converged}
  \STATE Compute norm, momenta, and kernel coefficients from \eqref{eq:gsvd-solution}
  \STATE Compute standard unconstrained gradient $\nabla E$ from \eqref{eq:standard-gradient-formula} and \eqref{eq:wpnorm-variation-particle}
  \STATE Form projected natural gradient $\widetilde{\nabla E}_n$ from \eqref{eq:operator-Y-definition} and \eqref{eq:discrete-path-metric}.
  \STATE Find $\varepsilon$ such that $E(V - \varepsilon \widetilde{\nabla E}_n) < E(V)$
  \STATE Update $V \gets V - \varepsilon \widetilde{\nabla E}_n$, update particles using \eqref{eq:particle-update}
\ENDWHILE
\STATE Output: Velocity field $V$
\end{algorithmic}

\caption{A simple gradient descent algorithm for optimization.}\label{alg:min-algorithm}
\end{algorithm}

Our stopping criterion for the optimization iteration is defined by monitoring the relative objective decrease at each step in tandem with the Euclidean vector 2-norm of the projected natural gradient (normalized with respect to the size of the vector). At convergence, the former is less than $10^{-8}$ and the latter is less than $10^{-6}$.

One necessary measure of convergence is constancy of the $W P$ norm along the path. For all the results we show, the relative variation measure $(\max_t \|\V{v}^t\|^2_{W P(\lmpt)})/(\min_t \|\V{v}^t\|^2_{W P(\lmpt)}) - 1$ is no larger than $10^{-3}$ and is usually $\mathcal{O}(10^{-6})$.

In all the tests we use 150 particles. The full set of particles is obtained by computing 50 uniformly distributed particle samples on the ranges of the two welding maps $\phi_0$ and $\phi_1$, and 50 on the domain of the welding maps, and then taking the union of the sets. This encourages inclusion of resolvable features from both the interior and exterior of the welds in the computation.

For the temporal discretization, we have found that using a high-order Gauss-Lobatto quadrature rule works very well for welding maps whose derivatives are not large. When the welding maps exhibit high frequency content (a result of e.g. protrusions in the shape), then a standard piecewise linear choice works well. In this latter case we use 150 equispaced points in the time variable.

\subsection{Zero-distance welds}
\begin{figure}
  \begin{tabular*}{\textwidth}{m{0.40\textwidth}c}
    \resizebox{0.40\textwidth}{!}{
    \renewcommand{\arraystretch}{1.3}
    \setlength{\tabcolsep}{12pt}
      \begin{tabular}{@{}rcccc@{}}
        \multicolumn{5}{c}{\vspace*{-350pt}}\\
        \multicolumn{5}{c}{Computed path energy between welds in Figure \ref{fig:welding-maps}.}\\\toprule
                 & $\phi_1$        & $\phi_2$       & $\phi_3$        & $\phi_4$ \\
        $\phi_1$ & ---             & ---            & ---             & ---      \\
        $\phi_2$ & \num{7.530e-12} & ---            & ---             & ---      \\
        $\phi_3$ & \num{1.249e-9}  & \num{1.729e-8} & ---             & ---      \\
        $\phi_4$ & \num{2.853e-9}  & \num{6.633e-8} & \num{2.069e-10} & ---      \\
      \end{tabular}
    \renewcommand{\arraystretch}{1}
    \setlength{\tabcolsep}{6pt}
    }
    &
    \resizebox{0.55\textwidth}{!}{
      \includegraphics{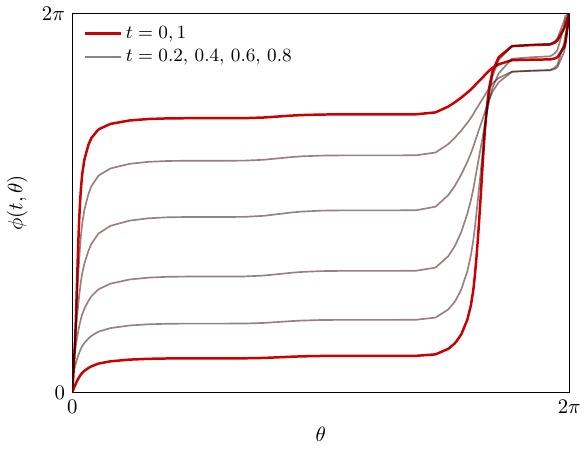}
    }
  \end{tabular*}
  \caption{Left: computed path energies between welding maps in the same equivalence
  class shown in Figure \ref{fig:welding-maps}. The near-zero value of the
  energies verifies that the computed paths are accurate geodesics. Right:
  evolution snapshots of $\phi(t,\theta)$ at $t=0, 0.2, 0.4, \ldots, 1$ where
  $\phi(0,\cdot) = \phi_3$ and $\phi(0,\cdot) = \phi_4$.}\label{fig:equivalence-class-energies}
\end{figure}

As a first test we compute the energies of the paths connecting welding maps in
the same equivalence class (which are expected to vanish). The welding maps are
the maps $\phi_i$ in Figure \ref{fig:welding-maps}. In Table \ref{fig:equivalence-class-energies} we show the computed energy at convergence. The shown near-vanishing energies confirms validity of the computation.

\subsection{Path length vs aspect ratio}
\begin{figure}
  \begin{center}
    \resizebox{1.0\textwidth}{!}{
      \includegraphics{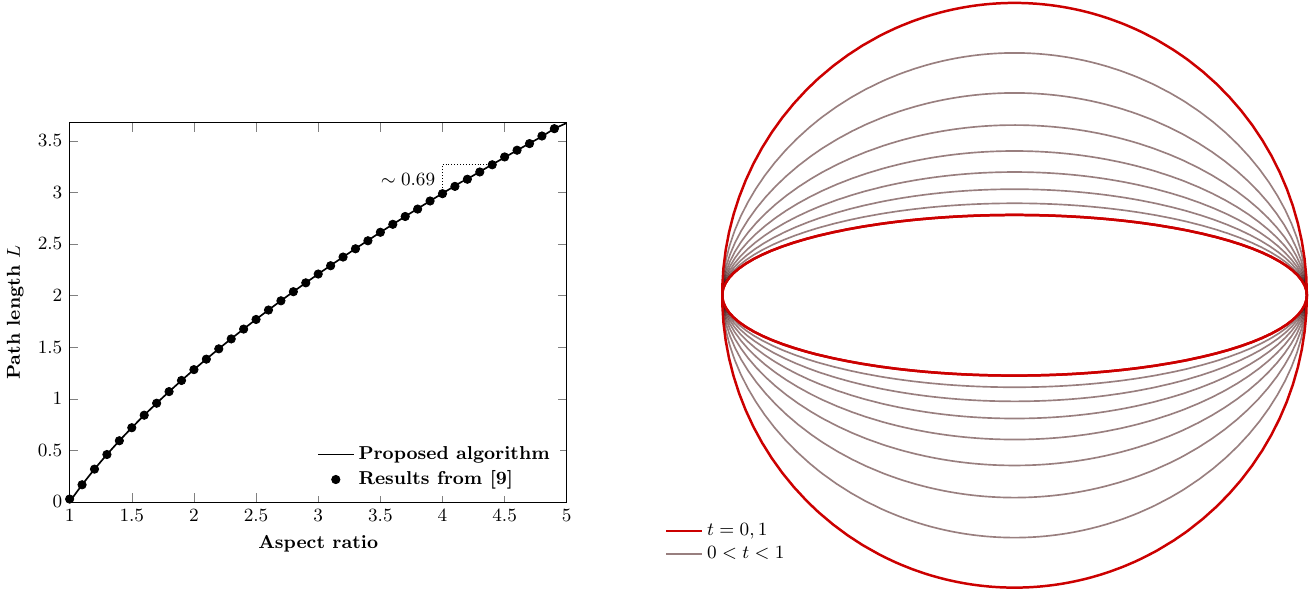}
    }
  \end{center}
  \caption{Left: computed geodesic energy for evolving an ellipse of given aspect ratio to a circle using the minimization algorithm (solid line) and using the shooting procedure from \cite{kushnarev-2012} (dashed line). (The lines overlap.) Right: evolution of the ellipse shown at equidistant points of the path parameter $s$.}
  \label{fig:ellipse-aspect-ratio}
\end{figure}

Let $\psi^{r}_t$ be a geodesic where $\psi^{r}_0$ is the welding map for an ellipse of aspect ratio $r$, and $\psi^{r}_1$ is the identity. In 
Figure \ref{fig:ellipse-aspect-ratio} we show results of simulations for path lengths $L[\psi^r]$ versus aspect ratio $r$ ranging from 1 to 5. The results suggest that the asymptotic relation is approximately linear, and these results match those obtained in \cite{kushnarev-2012} very well. 

\subsection{Hyperbolicity test}
\begin{figure}
  \begin{center}
    \resizebox{1.0\textwidth}{!}{
      \includegraphics{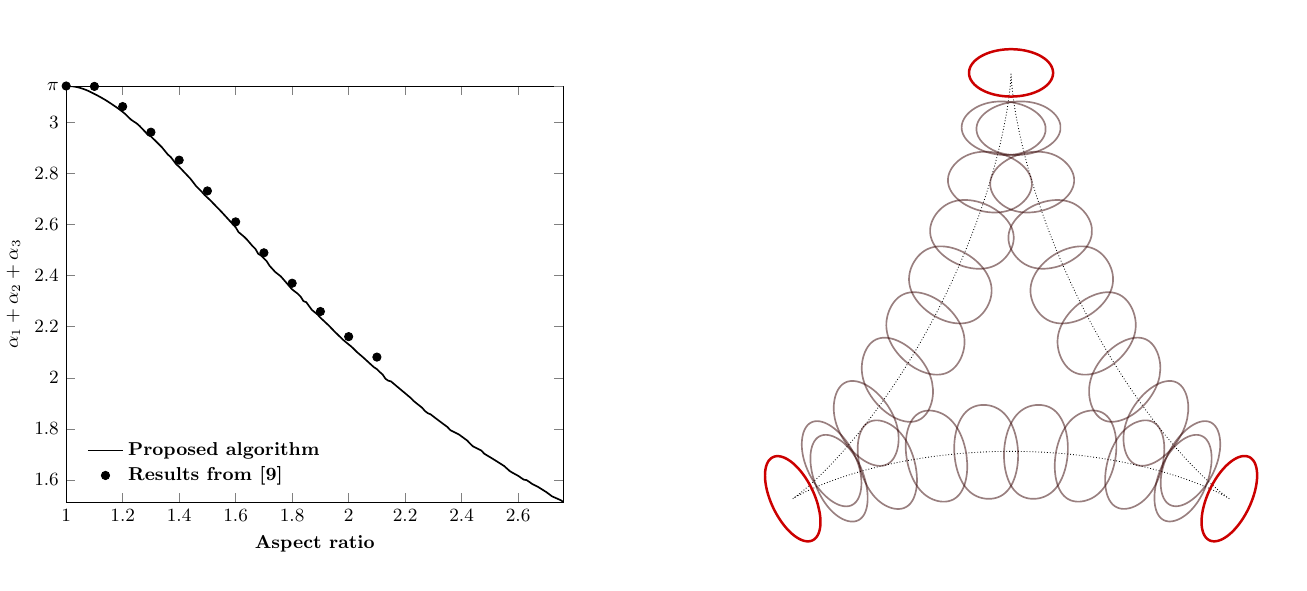}
    }
  \end{center}
  \caption{Left: Angle sum versus aspect ratio using the minimization algorithm and the results from \cite{kushnarev-2012}. Right: evolution along geodesics for a triplet of rotated ellipses. Each ellipse-like shape is the true shape along the geodesic, but the translational placement along a hyperbolic triangle is a cartoon-like representation of the flow on $T(1)$, suggesting the negative curvature of the space.}
  \label{fig:negative-curvature}
\end{figure}

We verify the negative curvature of Teichm\"uller space $\UTS$. Consider an ellipse of a fixed aspect ratio corresponding to the weld $\phi_1$. For rotations of the ellipse by angles $\pm 2 \pi/3$ we generate two more welding maps $\phi_2$ and $\phi_3$. We consider these three points in $\UTS$ and form a triangle using $\phi_i$ as the vertices and corresponding geodesic paths on $\UTS$ as the vertex connections. Let the $t=0$ velocity field that pushes $\phi_i$ to $\phi_j$ at $t=1$ be denoted $v_{i,j}$. Then the angle at vertex $\phi_i$ is given by
\begin{align*}
  \alpha_{i} = \frac{\langle v_{i,i\oplus 1}, v_{i, i\ominus 1} \rangle_{W P}}{\left\|v_{i,i\oplus 1}\right\|_{W P} \left\|v_{i,i\ominus 1}\right\|_{W P}},
\end{align*}
where $\oplus$ and $\ominus$ denote modular addition and subtraction, respectively, on the set $\{1, 2, 3\}$. Because all sectional curvatures of the $W P$ metric are negative, we expect that $\sum_{i=1}^3 \alpha_i \leq \pi$. We verify this fact in Figure \ref{fig:negative-curvature}. We again compare the results against those computed in \cite{kushnarev-2012} and obtain similar results.

\subsection{Shapes with corners}
\begin{figure}
  \begin{center}
    \resizebox{1.0\textwidth}{!}{
      \includegraphics{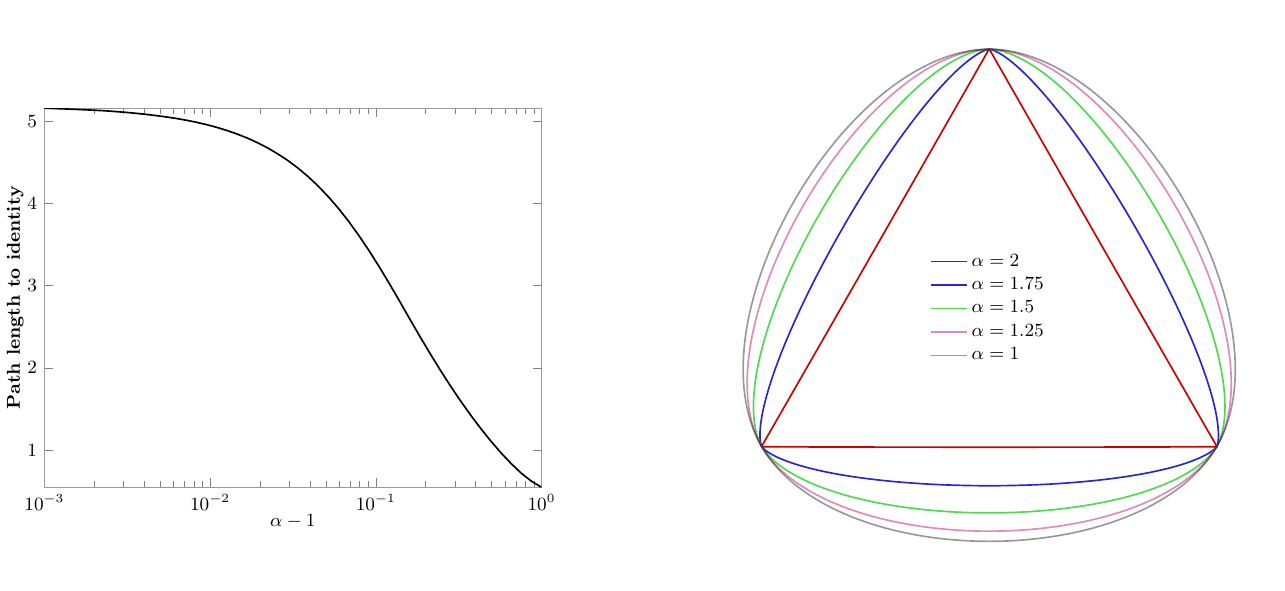}
    }
  \end{center}
  \caption{Left: $W P$ distance from a circle (the identity weld) to a triangle with smoothness parameter $\alpha$. Right: triangular shapes corresponding to values of the $\alpha$ parameter.}
  \label{fig:sobolev-smoothness}
\end{figure}

Non-smooth shapes do not result in diffeomorphic welds, and so the distance from any element of $T(1)$ to these non-smooth welds in the $W P$ metric is infinite. Our algorithm cannot compute a geodesic of infinite length, so we expect to run into limitations in the algorithm as we attempt to compute geodesics to shapes with sharp corners. We generate a triangular shape with rounded corners in the plane, and parameterize the smoothness of the corner by $\alpha \geq 1$. As shown on the right-hand side of Figure \ref{fig:sobolev-smoothness}, $\alpha=1$ corresponds to a sharp corner (a triangle in Euclidean space), and increasing $\alpha$ rounds the corners of the triangle. For many values of $\alpha$, we compute the geodesic from this rounded triangle to the identity and compile how the path length depends on $\alpha$. The results are shown on the left-hand side of Figure \ref{fig:sobolev-smoothness}. 

As we decrease $\alpha$ down to the critical value of $1$, we do observe a sharp growth in the path length. When $\alpha$ gets close to $1 + 10^{-3}$, our algorithm begins to reach its limit: at termination of the algorithm, the relative variation of the $W P$ norm along the path increases to about $0.1\%$. Results from computations with $\alpha -1 $ smaller than $10^{-3}$ no longer appear to be geodesics. (I.e., the $W P$ norm variation along the path becomes larger.) While this illustrates limitations of the algorithm, we are able to verify the expected increase in path length as $\alpha \downarrow 1$.

\subsection{Shapes from the MPEG-7 dataset}
\begin{figure}
  \begin{center}
    \resizebox{1.0\textwidth}{!}{
      \includegraphics{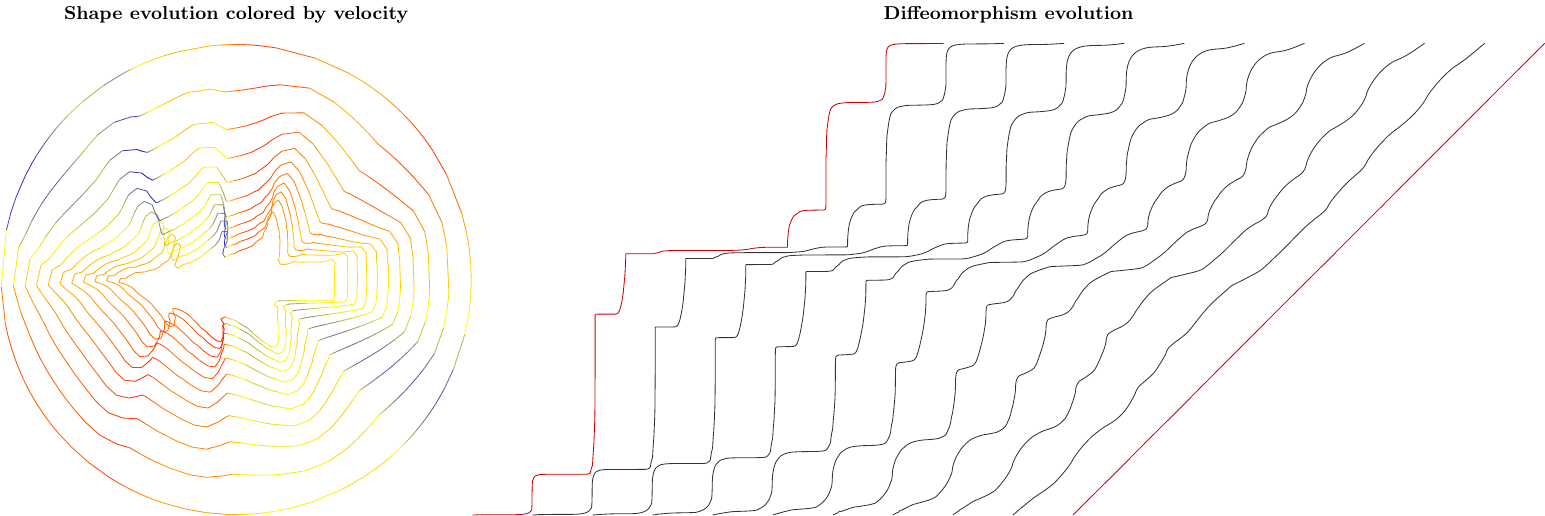}
    }
  \end{center}
  \caption{Left: shape evolution for {\sc mpeg-7} shape 1333, a rotated tree. Right: Evolution of the fingerprint along the path. The path length on $T(1)$ is $10.25$.} 
  \label{fig:mpeg7-1333}
\end{figure}

\begin{figure}
  \begin{center}
    \resizebox{1.0\textwidth}{!}{
      \includegraphics{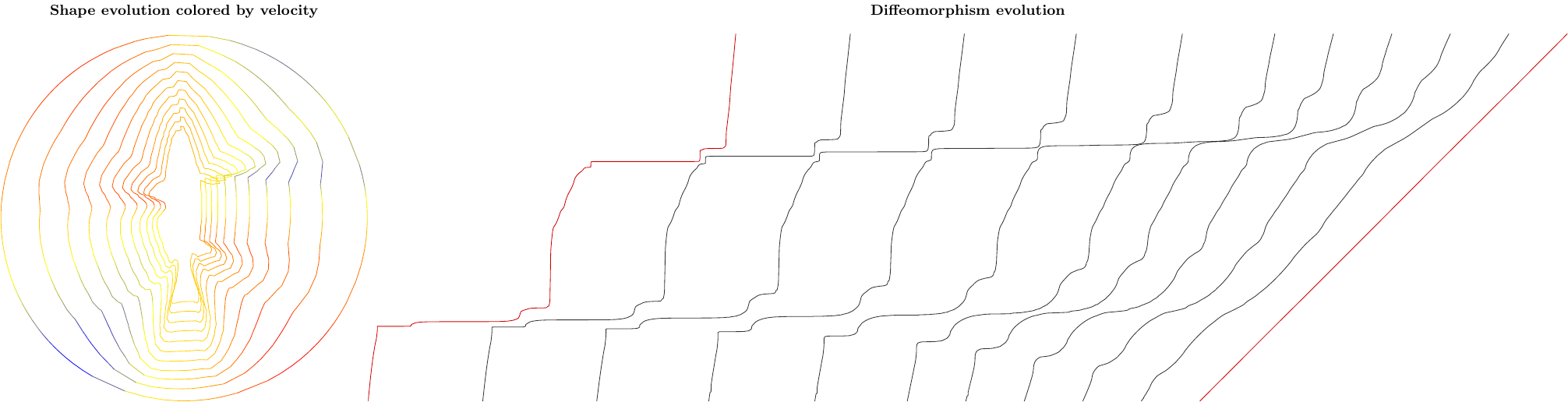}
    }
  \end{center}
  \caption{Left: shape evolution for {\sc mpeg-7} shape 785, a fish. Right: Evolution of the fingerprint along the path.. The path length on $T(1)$ is $8.889$.}
  \label{fig:mpeg7-785}
\end{figure}

\begin{figure}
  \begin{center}
    \resizebox{1.0\textwidth}{!}{
      \includegraphics{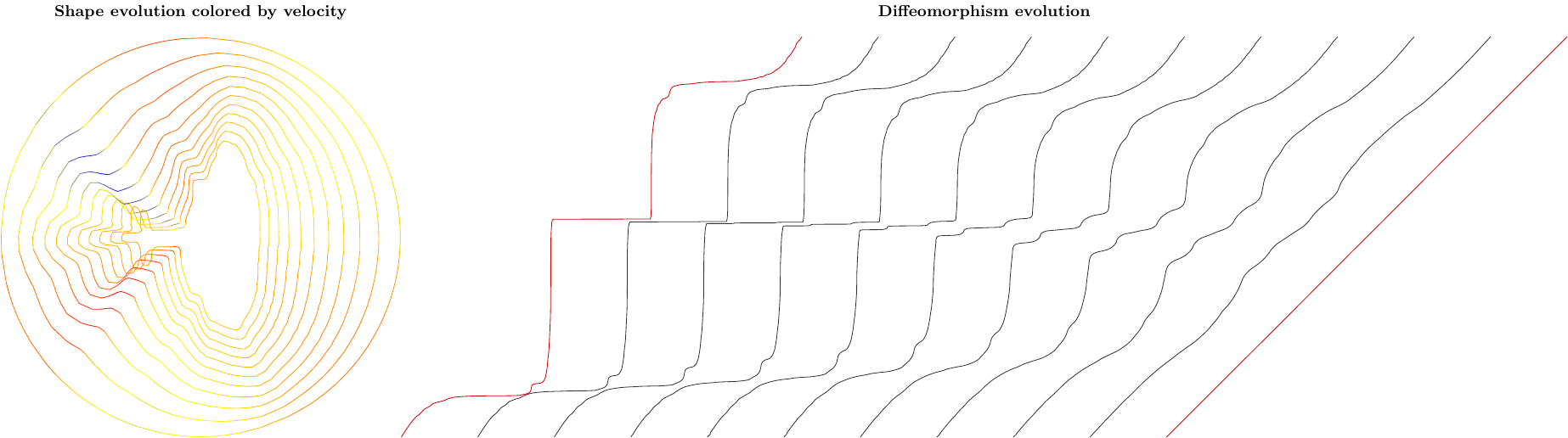}
    }
  \end{center}
  \caption{Left: shape evolution for {\sc mpeg-7} shape 875, a rotated miter. Right: Evolution of the fingerprint along the path. The path length on $T(1)$ is $7.486$.}
  \label{fig:mpeg7-875}
\end{figure}

We conclude our investigation with computation of geodesics from shapes in the MPEG-7 CE-shape-1 collection of planar shapes\cite{mpeg7}. Figures \ref{fig:mpeg7-1333}, \ref{fig:mpeg7-785}, and \ref{fig:mpeg7-875} show computed geodesics, and list the path length in each case. Also displayed are snapshots of the welding maps $\phi(t,\theta)$ at equidistant points in time.

We note that many shapes in this database have welding maps whose derivatives are singular to machine precision, or vanish to machine precision.\footnote{In this case the composite conformal maps of the weld are `crowded'; the max-min ratio of the derivative magnitude is large. A robust strategy for finite-precision computation in this situation is an open problem.}

When this happens, any algorithm implemented in finite precision for computing geodesics from this welding map will fail. We anticipate that any future work aimed at fixing the problem will first need to address the difficulty in computing conformal welds for these problematic shapes.

\section{Conclusion}
In this paper we have proposed and implemented a numerical solver for computing geodesics on the universal Teichm\"uller space $T(1)$ in the $W P$ metric. Given a planar shape, a conformal weld in the coset space $\UTS$ uniquely identifies this shape. We prescribe two given welding maps as the endpoints of the path in $\UTS$, and minimize the $W P$ energy of the path, subject to the endpoint constraints. This provides a constructive way to metrize the space of planar shapes, modulo rigid translation and scaling. 

We demonstrated the applicability of our method by solving for numerous geodesics; our results compare favorably with the results obtained by a shooting method in \cite{kushnarev-2012}. The main difficulty with our procedure stems from the inability to handle welds whose composite conformal maps contain large derivatives; this is in fact a limitation for any procedure of which we are aware. Our future work will be directed at addressing this problem.

Our optimization algorithm can, in effect, be used to minimize energy in any Riemannian metric, and our approach for satisfying endpoint constraints on manifolds can likewise be utilized in other applications.

\bibliographystyle{amsplain}
\bibliography{wp-bib}

\end{document}